\documentclass{myclass}
\usepackage[comma,numbers,square,sort&compress]{natbib}
\usepackage{epstopdf}
\usepackage{mathtools}
\usepackage{bm}
\usepackage{amsthm}
\usepackage[T1]{fontenc} 
\usepackage{amssymb}
\usepackage{xcolor}

\begin{document}
	
	\ArticleType{RESEARCH PAPER}
	\Year{2027}
	\Month{January}
	\Vol{70}
	\No{1}
	\DOI{}
	\ArtNo{}
	\ReceiveDate{}
	\ReviseDate{}
	\AcceptDate{}
	\OnlineDate{}
	\AuthorMark{}
	\AuthorCitation{}
	
	\title{On Observability Decomposition of Boolean Control Networks}{}
	
	\author[]{Mengjiao Liu}{250901008@njnu.edu.cn}
	\author[]{Jiandong Zhu*}{{zhujiandong@njnu.edu.cn}}

	
	\address[]{Key Laboratory for NSLSCS of Ministry of Education, School of Mathematical Sciences, \\ Nanjing Normal University, Nanjing 210023, China}
	
	\abstract{This paper studies the observability decomposition of Boolean networks (BNs) and Boolean control networks (BCNs) in the semi-tensor product framework. We characterize when a BN admits an observability decomposition in terms of the equivalence classes induced by indistinguishability. We prove that a decomposition exists if and only if all such equivalence classes have the same cardinality, and we explicitly construct the coordinate transformation together with the matrices of the observable and unobservable subsystems. The observable subsystem is shown to be the quotient graph obtained by contracting each class into a single vertex. For BNs whose equivalence classes have unequal cardinalities, we propose an extended observability decomposition by adding virtual states, and show that every finite BN admits an observability decomposition in an extended sense. These results are further generalized to BCNs, for which an algorithm and an illustrative example are provided.}

	\keywords{Boolean networks, Boolean control networks, indistinguishability classes, observability decomposition, extended decomposition}
	
	\maketitle

	\section{Introduction}
	Boolean networks (BNs) have attracted considerable attention in the natural sciences and in artificial systems. In a BN, each variable is binary, which makes BNs suitable for modeling systems with two-state components. The concept was originally proposed by Kauffman \cite{kauffman1969metabolic} to describe and characterize the dynamical behavior of gene regulatory networks, in which regulatory genes are either ``on'' or ``off''. A BN is an autonomous system: given the initial condition, its trajectory is uniquely determined.  When external Boolean control inputs are present, BNs are generalized to Boolean control networks (BCNs). Much research has been conducted on BCNs, including robustness \cite{aracena2009robustness} and synchronization \cite{hung2006chaos}. However, research on Boolean systems has long lacked a unified tool.
	
	In 2009, Cheng et al. proposed the semi-tensor product (STP) of matrices and established an algebraic framework for BNs and BCNs \cite{cheng2009input,cheng2009controllability,cheng2010linear,cheng2010state}. The core idea is to represent each of the \(2^n\) state combinations of an \(n\)-variable BN by a canonical basis vector in \(\mathbb{R}^{2^n}\). The logical mapping that governs the state update then becomes a \(2^n \times 2^n\) logical matrix. Consequently, a BN can be represented as a discrete-time linear system, and a BCN as a discrete-time bilinear system. The monograph \cite{cheng2010analysis} collects related results and treats controllability, observability, controller design, stability, disturbance decoupling, and other problems in this framework. Under this algebraic framework, logical problems are converted into algebraic ones and can be analyzed with algebraic state-space representation tools. For a fixed input, a BCN reduces to a BN; hence BCNs can be viewed as a family of BNs and studied by graph-theoretic methods \cite{laschov2012controllability,zhao2010input}. Subsequent work has addressed controllability \cite{li2011controllability,zhang2019set}, stabilization \cite{cheng2011stability,lu2019pinning}, disturbance decoupling \cite{cheng2010disturbance,li2019disturbance}, output regulation \cite{li2016output,zhang2019output}, switched models \cite{li2015controllability,li2017set,gao2020infinite}, and time-delay models \cite{zhang2013controllability,zhang2020event,mu2024identification}.
	
	Observability decomposition is a classical topic in control theory. For BCNs, observability and its criteria have been studied extensively \cite{cheng2009controllability,fornasini2012observability}, whereas observability decomposition has received comparatively less attention \cite{cheng2010realization,zou2017graph,li2022observability}. The work in \cite{cheng2010realization} initiated the study of this problem for BCNs by introducing regular subspaces: the decomposition is possible when the maximal unobservable subspace is regular and satisfies the regularity condition. This condition, however, is often difficult to verify. In \cite{zou2017graph}, the output decomposition problem was addressed by a vertex partition method, which identifies the largest unobservable part and separates it from the rest of the network. This removes the regularity assumption, but the resulting subsystem is not guaranteed to be observable. Building on the vertex partition method in \cite{zou2017graph}, the work in \cite{li2022observability} imposed additional conditions to obtain an observable subsystem; however, the resulting criterion for the realizability of observability decomposition involves several conditions. These observations motivate us to construct the decomposition directly from the indistinguishability equivalence classes.
	
	Following this line of research, we develop the decomposition directly from the indistinguishability classes. We show that a BN admits an observability decomposition if and only if all equivalence classes induced by indistinguishability have the same cardinality, and we explicitly construct the transformation matrix and the structure matrices for the decomposed system. The observable subsystem is exactly the quotient graph obtained by contracting each equivalence class into a single vertex. For BNs with unequal class cardinalities, we propose an extended decomposition by adding virtual states, and prove that every finite BN can be embedded into an extended logical system that admits an observability decomposition. We then generalize the results to BCNs and provide an algorithm and an illustrative example. In contrast to \cite{cheng2010realization,zou2017graph,li2022observability}, the proposed approach does not require the maximal unobservable subspace to be regular, guarantees that the separated subsystem is observable, and derives the coordinate transformation directly from the classes, which makes the quotient-graph structure explicit.
	
	The rest of this paper is organized as follows. Section 2 introduces the notations and the semi-tensor product. Section 3 establishes the criterion for observability decomposition of BNs via indistinguishability classes, and proposes the extended decomposition for BNs with classes of unequal cardinalities. Section 4 generalizes the results to BCNs, and provides an algorithm and an illustrative example. Section 5 concludes the paper.

	\section{Preliminaries}
	
	Throughout this paper, we use the following notations:
	\begin{itemize}
		\item \(|S|\): the cardinality of a set \(S\);
		\item \([a,b] := \{a, a+1, \ldots, b\}\), where \(a \le b\) are integers;
		\item \(Col_i(A)\): the \(i\)-th column of a matrix \(A\); \(Col(A)\) denotes the set of all columns of \(A\);
		\item \(A^T\): the transpose of a matrix \(A\);
		\item \(\mathbb{R}^n\): the set of \(n\)-dimensional real column vectors; \(\mathbb{R}^{m \times n}\) denotes the set of all \(m \times n\) real matrices;
		\item \(\mathcal{B} := \{0,1\}\): the Boolean domain; \(\mathcal{B}^n\) denotes the set of \(n\)-dimensional Boolean vectors;
		\item \(\bm{1}_k\): the \(k\)-dimensional column vector with all entries equal to \(1\);
		\item \(\bm{I}_k\): the \(k \times k\) identity matrix;
		\item \(\delta_k^i\): the \(i\)-th column of \(\bm{I}_k\), called a \(k\)-dimensional logical vector; when the dimension is clear from the context, it is simply denoted by \(i\);
		\item \(\delta_n[i_1\ i_2\ \cdots\ i_r] := [\delta_n^{i_1}\ \delta_n^{i_2}\ \cdots\ \delta_n^{i_r}]\), where \(i_j \in [1,n]\) for \(j=1,\ldots,r\);
		\item \(\Delta_k := \{\delta_k^i \mid i=1,2,\ldots,k\}\);
		\item \(\mathcal{L}_{m \times n}\): the set of all \(m \times n\) logical matrices, i.e., matrices all of whose columns are logical vectors;
		\item \(\otimes\): the Kronecker product;
		\item \(\mathcal{N}(\cdot)\): the out-neighbor set of an element or a set; \(\mathcal{N}_s(\cdot)\) denotes the out-neighbor set in the subgraph \(\mathcal{G}_s\).
	\end{itemize}
	
	A permutation matrix is a matrix \(T \in \mathcal{L}_{n \times n}\) with \(Col(T) = \Delta_n\); such a matrix is nonsingular and satisfies \(T^{-1} = T^T\).
	
	Given a matrix \(L \in \mathcal{L}_{n \times n}\), we associate with it a directed graph \(\mathcal{G}(L)\). There is an arc \((i,j)\) from \(i\) to \(j\) if and only if the \((j,i)\)-th entry of \(L\) is \(1\), i.e., \(\delta_n^j = L\delta_n^i\).
	
	There is a one-to-one correspondence between the Boolean variable \(X \in \mathcal{B}\) and the logical variable \(x \in \Delta_2\):
	\begin{equation}
		x = \begin{bmatrix} X \\ 1-X \end{bmatrix}.
		\label{oto}
	\end{equation}
	
	We introduce the semi-tensor product (STP) of matrices, denoted by \(\ltimes\), as follows \cite{cheng2010linear}: given \(A \in \mathbb{R}^{m \times n}\) and \(B \in \mathbb{R}^{p \times q}\),
	\begin{equation*}
		A \ltimes B := (A \otimes \bm{I}_{t/n})(B \otimes \bm{I}_{t/p}),
	\end{equation*}
	where \(t\) is the least common multiple of \(n\) and \(p\). When \(n=p\), \(\ltimes\) reduces to the conventional matrix multiplication. For more details on the STP, see \cite{cheng2010analysis}. In the sequel, we write \(\ltimes\) simply as juxtaposition when no confusion arises.
	
	With the aid of the STP, there is a one-to-one correspondence, similar to \eqref{oto}, between an \(n\)-dimensional Boolean vector \(X=[X_1,X_2,\ldots,X_n]^T \in \mathcal{B}^n\) and an \(N\)-dimensional logical vector \(x \in \Delta_N\) (\(N=2^n\)):
	\begin{equation*}
		x := \begin{bmatrix} X_1 \\ 1-X_1 \end{bmatrix} \ltimes \begin{bmatrix} X_2 \\ 1-X_2 \end{bmatrix} \ltimes \cdots \ltimes \begin{bmatrix} X_n \\ 1-X_n \end{bmatrix}.
	\end{equation*}
	
	\section{Observability Decomposition of BNs}
	
	A Boolean network is described by
	\begin{equation}
		\label{eq:BNLE}
		\begin{aligned}
			X(t+1) &= f(X(t)), \\
			Y(t) &= h(X(t)),
		\end{aligned}
	\end{equation}
	where \(X\in\mathcal{B}^n\) is the state variable, \(Y\in\mathcal{B}^p\) is the output variable, and \(f:\mathcal{B}^n\to\mathcal{B}^n\) and \(h:\mathcal{B}^n\to\mathcal{B}^p\) are logical functions.
	
	Let \(x\in\Delta_N\) and \(y\in\Delta_P\) represent \(X\) and \(Y\), respectively, where \(N=2^n\) and \(P=2^p\). Then BN \eqref{eq:BNLE} can be expressed in algebraic form as \cite{cheng2010analysis}:
	\begin{equation}
		\label{eq:BNAE}
		\begin{aligned}
			x(t+1) &= L \ltimes x(t) = L x(t), \\
			y(t) &= H \ltimes x(t) = H x(t),
		\end{aligned}
	\end{equation}
	where \(L\in\mathcal{L}_{N\times N}\) and \(H\in\mathcal{L}_{P\times N}\).
	
	\begin{definition}[{\cite{fornasini2012observability}}]
		Given a Boolean network \eqref{eq:BNAE},
		\begin{itemize}
			\item Two states \(x=\delta_N^i\) and \(\bar{x}=\delta_N^j\) are said to be indistinguishable if the output trajectories initialized at \(x\) and \(\bar{x}\), respectively, coincide for all \(t\ge 0\); otherwise, they are said to be distinguishable;
			\item Moreover, a BN is said to be observable if any two distinct states are distinguishable.
		\end{itemize}
	\end{definition}
	
	To analyze the observability decomposition problem, we introduce an equivalence relation on \(\Delta_N\): two states are related by \(\sim\) if they are indistinguishable. Since \(\sim\) is an equivalence relation, it partitions \(\Delta_N\) into equivalence classes; let \(\mathcal{C}_\sim\) denote the set of classes. We also say that \(x\) and \(\bar{x}\) are \(k\)-step indistinguishable, denoted by \(x\sim_k\bar{x}\), if the output trajectories initialized at \(x\) and \(\bar{x}\), respectively, coincide for all \(t\in[0,k-1]\). The relation \(\sim_k\) is again an equivalence relation, and we let \(\mathcal{C}_{\sim,k}\) denote the set of equivalence classes. Since the partition induced by \(\sim_{k+1}\) refines that induced by \(\sim_k\) and \(\Delta_N\) is finite, there exists \(k\in\mathbb{Z}^+\) such that \(\mathcal{C}_{\sim,k}=\mathcal{C}_{\sim,k+1}\); for such \(k\), one has \(\mathcal{C}_{\sim,k}=\mathcal{C}_{\sim,k+1}=\cdots=\mathcal{C}_\sim\). For more details on the equivalence classes induced by indistinguishability, see \cite{fornasini2012observability}.
	
	\begin{figure}[!t]
		\centering
		\begin{minipage}[c]{0.48\textwidth}
			\centering
			\includegraphics[width=3.91cm,height=3.6cm]{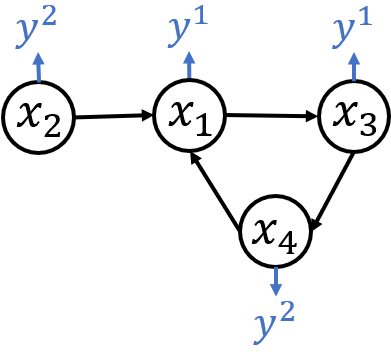}
		\end{minipage}
		\hspace{0.02\textwidth}
		\begin{minipage}[c]{0.48\textwidth}
			\centering
			\includegraphics[width=3.92cm,height=2.00cm]{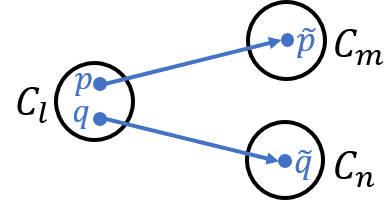}
		\end{minipage}\\[3mm]
		\begin{minipage}[t]{0.48\textwidth}
			\centering
			\caption{State transition graph with outputs of Example~\ref{eg1}.}
			\label{fig1}
		\end{minipage}
		\hspace{0.02\textwidth}
		\begin{minipage}[t]{0.48\textwidth}
			\centering
			\caption{Illustration of the proof process for Proposition~\ref{prop1}.}
			\label{fig2}
		\end{minipage}
	\end{figure}
	
	\begin{example}
		\label{eg1}
		Consider the BN described by the directed graph in Fig.~\ref{fig1}, where the circles represent states, the black arrows represent transitions, and the blue arrows represent outputs. Let \(x_i\) denote the state \(\delta_4^i\), \(i=1,\ldots,4\). From Fig.~\ref{fig1}, the one-step equivalence classes are
		\[
		\mathcal{C}_{\sim,1}=\bigl\{\{x_2,x_4\},\{x_1,x_3\}\bigr\}.
		\]
		The two-step equivalence classes are
		\[
		\mathcal{C}_{\sim,2}=\bigl\{\{x_2,x_4\},\{x_1\},\{x_3\}\bigr\},
		\]
		and the three-step equivalence classes are consistent with the two-step ones, i.e., \(\mathcal{C}_{\sim,3}=\mathcal{C}_{\sim,2}\). Hence the equivalence classes induced by \(\sim\) are
		\[
		\mathcal{C}_\sim=\bigl\{\{x_2,x_4\},\{x_1\},\{x_3\}\bigr\}.
		\]
	\end{example}
	
	\begin{proposition}
		\label{prop1}
		Given BN \eqref{eq:BNAE}, let \(C_1,C_2,\ldots,C_r\) be the equivalence classes induced by the indistinguishability relation. Then, for any \(C_i\), there exists an index \(k\) such that
		\[
		\mathcal{N}(C_i)\subseteq C_k.
		\]
	\end{proposition}
	
	\begin{proof}
		Assume, for contradiction, that there exists some \(C_i\) for which there exist \(C_m\) and \(C_n\) with \(m\neq n\) such that
		\[
		\mathcal{N}(C_i)\cap C_m\neq\emptyset,\quad \mathcal{N}(C_i)\cap C_n\neq\emptyset.
		\]
		Choose \(p,q\in C_i\) and \(\tilde{p} \in \mathcal{N}(p)\cap C_m\), \(\tilde{q}\in \mathcal{N}(q)\cap C_n\). Since \(p\) and \(q\) belong to the same class \(C_i\), they are indistinguishable, so the output trajectories from \(p\) and \(q\) coincide for all \(t\ge 0\). At \(t=1\), these outputs are \(H \tilde{p}\) and \(H \tilde{q}\), hence \(H \tilde{p} = H \tilde{q}\). Moreover, for \(t\ge 1\), the outputs from \(p\) and \(q\) are determined by the trajectories from \(\tilde{p}\) and \(\tilde{q}\), respectively. Therefore, \(\tilde{p}\) and \(\tilde{q}\) are indistinguishable. But \(\tilde{p}\in C_m\) and \(\tilde{q}\in C_n\) with \(m\neq n\), so they are distinguishable, a contradiction. This is illustrated in Fig.~\ref{fig2}.
		\hfill\(\square\)
	\end{proof}
	
	Before defining observability decomposition for BN \eqref{eq:BNAE}, we recall the corresponding notion for linear systems. Consider
	\begin{equation}
		\begin{cases}
			\dot{x} = Ax, \\
			y = Cx,
		\end{cases}
		\quad
		\xleftrightarrow[\scriptstyle x=T^{-1}z]{\scriptstyle z=Tx}
		\quad
		\begin{cases}
			\dot{z}_1 = A_{11}z_1, \\
			\dot{z}_2 = A_{21}z_1 + A_{22}z_2, \\
			y = C_1 z_1,
		\end{cases}
		\label{eq:linear_decomp}
	\end{equation}
	where the pair \((C_1,A_{11})\) is observable. Here \(x\in\mathbb{R}^n\) and \(y\in\mathbb{R}^p\) are the state and output of the original system, respectively, \(z_1\in\mathbb{R}^{n_1}\), \(z_2\in\mathbb{R}^{n_2}\), and \(T=[\beta_1\ \cdots\ \beta_{n_1}\mid \beta_{n_1+1}\ \cdots\ \beta_n]\) is an invertible matrix such that the vectors \(\beta_{n_1+1},\ldots,\beta_n\) span the unobservable subspace.
	
	\begin{definition}[{\cite{li2022observability}}]
		\label{def2}
		Consider BN \eqref{eq:BNAE}. An observability decomposition is said to be realizable if there exists a coordinate transformation \(z=Tx\) with a permutation matrix \(T\) such that, in the \(z\)-coordinate frame, the BN is transformed into
		\begin{equation}
			\begin{aligned}
				z^{[1]}(t+1) &= G_1 z^{[1]}(t), \\
				z^{[2]}(t+1) &= G_2 z(t), \\
				y(t) &= Q z^{[1]}(t),
			\end{aligned}
			\label{eq:BN_decomp}
		\end{equation}
		where the \(z^{[1]}\)-subsystem is observable, \(z^{[1]}\in\Delta_{N_1}\), \(z^{[2]}\in\Delta_{N_2}\), \(y\in\Delta_P\), \(z=z^{[1]}\ltimes z^{[2]}\), and \(N_1=2^{n_1}\), \(N_2=2^{n_2}\), \(P=2^p\), \(N_1N_2=N\). The matrices satisfy \(G_1\in\mathcal{L}_{N_1\times N_1}\), \(G_2\in\mathcal{L}_{N_2\times N}\), and \(Q\in\mathcal{L}_{P\times N_1}\).
	\end{definition}
	
	\begin{remark}
		When \(n_1=n\), BN \eqref{eq:BNAE} is observable; when \(n_2=n\), all states in BN \eqref{eq:BNAE} are indistinguishable.
	\end{remark}
	
	\begin{remark}
		The decomposition \eqref{eq:BN_decomp} has the same structure as \eqref{eq:linear_decomp}: the observable part is determined only by \(z^{[1]}\), while the unobservable part is determined by the whole state \(z\).
	\end{remark}
	
	We use \((L,H)\) and \((G_1,G_2,Q)\) to represent BN \eqref{eq:BNAE} and \eqref{eq:BN_decomp}, respectively. The coordinate transformation \(z=Tx\) gives the following lemma.
	
	\begin{lemma}[{\cite{zou2017graph}}]
		If BN \eqref{eq:BNAE} admits an observability decomposition, then \((L,H)\) and \((G_1,G_2,Q)\) satisfy
		\begin{align}
			&(\bm{I}_{N_1}\otimes\bm{1}_{N_2}^T)TLT^T
			=G_1(\bm{I}_{N_1}\otimes\bm{1}_{N_2}^T), \nonumber \\
			&(\bm{1}_{N_1}^T\otimes\bm{I}_{N_2})TLT^T=G_2, \nonumber \\
			&HT^T=Q(\bm{I}_{N_1}\otimes\bm{1}_{N_2}^T). \nonumber
		\end{align}
	\end{lemma}
	
	In the STP framework, the individual logical components in \eqref{eq:BNLE} are combined into the overall variable \(x\) in \eqref{eq:BNAE}. The following lemma shows that each component can also be extracted from \(x\).
	
	\begin{lemma}[{\cite{cheng2010analysis}}]
		\label{lemma2}
		\begin{itemize}
			\item Let \(x_i\in\Delta_2\) for \(i=1,\ldots,n\) and \(x=\ltimes_{i=1}^n x_i\), then
			\begin{equation*}
				x_i=(\bm{1}_{2^{i-1}}^T\otimes\bm{I}_2\otimes\bm{1}_{2^{n-i}}^T)x;
			\end{equation*}
			\item Similarly, the partial state from the \((i+1)\)-th to the \((i+k)\)-th block can also be obtained from \(x\):
			\begin{equation*}
				\ltimes_{j=i+1}^{i+k}x_j
				=(\bm{1}_{2^i}^T\otimes\bm{I}_{2^k}\otimes\bm{1}_{2^{n-i-k}}^T)x.
			\end{equation*}
		\end{itemize}
	\end{lemma}
	
	\begin{theorem}
		\label{thm1}
		BN \eqref{eq:BNAE} admits an observability decomposition if and only if all equivalence classes induced by the indistinguishability relation have the same cardinality.
	\end{theorem}
	
	\begin{proof}
		\emph{[Necessity]}
		Suppose BN~\eqref{eq:BNAE} admits an observability decomposition~\eqref{eq:BN_decomp}. For \(l=1,\ldots,N_1\), define
		\[
		C_l=\bigl\{\delta_N^a\in\Delta_N\mid
		\delta_N^a=T^T(\delta_{N_1}^l\ltimes\delta_{N_2}^b),\ b\in[1,N_2]\bigr\}.
		\]
		Since \(T\) is a permutation matrix and
		\[
		\{\delta_{N_1}^l\ltimes\delta_{N_2}^b\mid l\in[1,N_1],\ b\in[1,N_2]\}=\Delta_N,
		\]
		the sets \(C_1,\ldots,C_{N_1}\) form a partition of \(\Delta_N\), and \(|C_l|=N_2\) for every \(l\).
		
		We show that the sets \(C_l\) are exactly the equivalence classes induced by indistinguishability.
		
		First, let \(x_1,x_2\in C_l\). By the definition of \(C_l\), their \(z^{[1]}\)-coordinates are both \(\delta_{N_1}^l\). In the decomposition \eqref{eq:BN_decomp}, the \(z^{[1]}\)-subsystem evolves as \(z^{[1]}(t+1)=G_1z^{[1]}(t)\), and the output is \(y(t)=Qz^{[1]}(t)\). Hence the two states produce the same \(z^{[1]}\)-trajectory and the same output trajectory, so they are indistinguishable.
		
		Second, let \(x_1\in C_l\) and \(x_2\in C_k\) with \(l\neq k\). Their \(z^{[1]}\)-coordinates are \(\delta_{N_1}^l\) and \(\delta_{N_1}^k\), respectively. Since the \(z^{[1]}\)-subsystem is observable, \(\delta_{N_1}^l\) and \(\delta_{N_1}^k\) are distinguishable; that is, the corresponding output trajectories differ at some time instant. Therefore, \(x_1\) and \(x_2\) are distinguishable.
		
		Thus \(C_1,\ldots,C_{N_1}\) are the equivalence classes induced by indistinguishability, and all of them have cardinality \(N_2\).
		
		\emph{[Sufficiency]}
		Suppose the equivalence classes induced by indistinguishability are \(C_1,\ldots,C_{N_1}\), and \(|C_l|=N_2\) for \(l=1,\ldots,N_1\), where \(N_1N_2=N\). We construct \(T,G_1,G_2,Q\) satisfying \eqref{eq:BN_decomp}.
		
		\emph{Construction of \(T\).}
		For each \(l\), write
		\[
		C_l = \{ \delta_N^{a_1^l}, \delta_N^{a_2^l}, \cdots, \delta_N^{a_{N_2}^l} \}.
		\]
		Define
		\begin{equation}
			\label{eq:t}
			Col_{a_i^l}(T)=\delta_{N_1}^l\ltimes\delta_{N_2}^i, \quad l=1,\ldots,N_1,\ i=1,\ldots,N_2.
		\end{equation}
		Since \(C_1,\ldots,C_{N_1}\) form a partition of \(\Delta_N\), every column index appears exactly once, so \(T\) is a permutation matrix. Moreover,
		\begin{equation}
			T\delta_N^{a_i^l}=\delta_{N_1}^l\ltimes\delta_{N_2}^i.
			\label{t}
		\end{equation}
		Let \(z=Tx\). By Lemma~\ref{lemma2},
		\begin{equation}
			z^{[1]}(t)=(\bm{I}_{N_1}\otimes\bm{1}_{N_2}^T)Tx(t),
			\label{eq:z1}
		\end{equation}
		\begin{equation}
			z^{[2]}(t)=(\bm{1}_{N_1}^T\otimes\bm{I}_{N_2})Tx(t).
			\label{eq:z2}
		\end{equation}
		
		\emph{Construction of \(G_1\).}
		For each \(l\), by Proposition~\ref{prop1}, there exists an index \(k\) such that \(\mathcal{N}(C_l)\subseteq C_k\). Define
		\begin{equation}
			Col_l(G_1)=\delta_{N_1}^k,\quad l=1,\ldots,N_1.
			\label{eq:g1}
		\end{equation}
		This is independent of the representative chosen in \(C_l\), because the image of the whole class \(C_l\) under \(L\) lies in the single class \(C_k\). By \eqref{t} and \eqref{eq:z1}, if \(x(t)=\delta_N^{a_i^l}\), then \(z^{[1]}(t)=\delta_{N_1}^l\) and \(x(t+1)=L\delta_N^{a_i^l}\in C_k\), by \eqref{t}, \(z^{[1]}(t+1)=\delta_{N_1}^k=G_1\delta_{N_1}^l\). Hence \(z^{[1]}(t+1)=G_1z^{[1]}(t)\).
		
		\emph{Construction of \(Q\).}
		For each \(l\), all states in \(C_l\) have the same output. Let
		\begin{equation}
			Col_l(Q)=H\delta_N^{a_1^l},\quad l=1,\ldots,N_1.
			\label{eq:q}
		\end{equation}
		This definition is independent of the choice of the representative \(a_1^l\in C_l\). If \(x(t)=\delta_N^{a_i^l}\), then \(z^{[1]}(t)=\delta_{N_1}^l\) and
		\[
		y(t)=H\delta_N^{a_i^l}=H\delta_N^{a_1^l}=Q\delta_{N_1}^l=Qz^{[1]}(t).
		\]
		
		\emph{Construction of \(G_2\).}
		Let
		\begin{equation}
			G_2=(\bm{1}_{N_1}^T\otimes\bm{I}_{N_2})TLT^T.
			\label{eq:g2}
		\end{equation}
		Then
		\[
		z^{[2]}(t+1)
		\stackrel{\text{\eqref{eq:z2}}}{=}(\bm{1}_{N_1}^T\otimes\bm{I}_{N_2})Tx(t+1)
		=(\bm{1}_{N_1}^T\otimes\bm{I}_{N_2})TLx(t)
		=(\bm{1}_{N_1}^T\otimes\bm{I}_{N_2})TLT^T Tx(t)
		=G_2z(t).
		\]
		
		It remains to show that the \(z^{[1]}\)-subsystem is observable. Let \(l\neq k\). Suppose, for contradiction, that \(\delta_{N_1}^l\) and \(\delta_{N_1}^k\) are indistinguishable in the \(z^{[1]}\)-subsystem. Then for all \(t\ge 0\),
		\[
		Q(G_1)^t\delta_{N_1}^l=Q(G_1)^t\delta_{N_1}^k.
		\]
		Take any \(x_l\in C_l\) and \(x_k\in C_k\). By the construction of \(T\), their \(z^{[1]}\)-coordinates are \(\delta_{N_1}^l\) and \(\delta_{N_1}^k\), respectively. Hence their output trajectories in the original BN are
		\[
		y_l(t)=Q(G_1)^t\delta_{N_1}^l,\quad
		y_k(t)=Q(G_1)^t\delta_{N_1}^k,
		\]
		which coincide for all \(t\ge 0\). This contradicts the fact that \(x_l\) and \(x_k\) belong to different equivalence classes, hence are distinguishable. Therefore, any two distinct states in the \(z^{[1]}\)-subsystem are distinguishable, and the \(z^{[1]}\)-subsystem is observable.
		
		Therefore, the constructed \(T,G_1,G_2,Q\) satisfy \eqref{eq:BN_decomp} and the \(z^{[1]}\)-subsystem is observable. Hence BN \eqref{eq:BNAE} admits an observability decomposition.
		\hfill\(\square\)
	\end{proof}
	
	Below, we perform observability decomposition on a simple BN. This example illustrates Theorem~\ref{thm1} and its proof.
	
	\begin{example}
		\label{eg2}
		Consider the BN
		\begin{equation}
			\label{eq:eg2}
			\begin{aligned}
				x(t+1) &= \delta_4[2\ 3\ 4\ 1]x(t), \\
				y(t) &= \delta_2[1\ 2\ 1\ 2]x(t).
			\end{aligned}
		\end{equation}
		Its state transition graph with outputs is given by the black and blue parts in Fig.~\ref{fig3}~(a), where \(i\) represents \(\delta_4^i\), \(i=1,2,3,4\), and \(y_j\) represents \(\delta_2^j\), \(j=1,2\).
		
		\begin{figure}[!t]
			\centering
			\begin{minipage}[c]{0.48\textwidth}
				\centering
				\includegraphics[width=3.03cm,height=3.86cm]{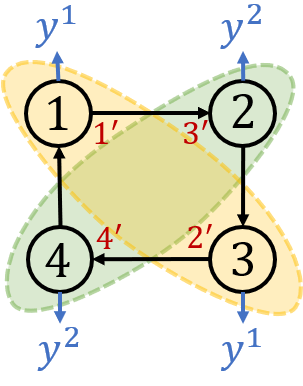}
			\end{minipage}
			\hspace{0.02\textwidth}
			\begin{minipage}[c]{0.48\textwidth}
				\centering
				\includegraphics[width=2.82cm,height=1.54cm]{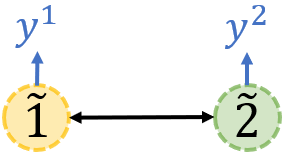}
			\end{minipage}
			\caption{(a) State transition graph of \eqref{eq:eg2}; (b) Quotient graph obtained by merging indistinguishable states.}
			\label{fig3}
		\end{figure}
		
		From Fig.~\ref{fig3}~(a), the equivalence classes induced by indistinguishability are \(\{1,3\}\) and \(\{2,4\}\), as shown in the yellow and green shaded regions in Fig.~\ref{fig3}~(a). Contracting \(\{1,3\}\) into a point \(\tilde{1}\) and \(\{2,4\}\) into a point \(\tilde{2}\), we obtain the quotient graph in Fig.~\ref{fig3}~(b). Let \(\tilde{1}\) and \(\tilde{2}\) be represented by \(\bigl[\begin{smallmatrix}1\\0\end{smallmatrix}\bigr]\) and \(\bigl[\begin{smallmatrix}0\\1\end{smallmatrix}\bigr]\), respectively. Then the state equation of the contracted graph is
		\begin{equation*}
			\begin{aligned}
				z^{[1]}(t+1) &= \delta_2[2\ 1]z^{[1]}(t),\\
				y(t) &= \delta_2[1\ 2]z^{[1]}(t).
			\end{aligned}
		\end{equation*}
		
		To obtain the decomposition \eqref{eq:BN_decomp}, we also need the evolution of \(z^{[2]}\). For a fixed value of \(z^{[1]}\), the set of semi-tensor products of this value with all values of \(z^{[2]}\) should cover exactly one equivalence class. Thus we assign \(z^{[2]}=\delta_2^1\) to states 1 and 2, and \(z^{[2]}=\delta_2^2\) to states 3 and 4. The coordinate correspondence is
		\begin{equation}
			\label{eq:xz1}
			\renewcommand{\arraystretch}{1.5}
			\setlength{\arraycolsep}{5pt}
			\begin{array}{ccccc}
				x & 1 & 2 & 3 & 4\\
				z^{[1]} & \bigl[\begin{smallmatrix}1\\0\end{smallmatrix}\bigr] & \bigl[\begin{smallmatrix}0\\1\end{smallmatrix}\bigr] & \bigl[\begin{smallmatrix}1\\0\end{smallmatrix}\bigr] & \bigl[\begin{smallmatrix}0\\1\end{smallmatrix}\bigr]\\
				z^{[2]} & \delta_2^1 & \delta_2^1 & \delta_2^2 & \delta_2^2\\
				z & 1' & 3' & 2' & 4'\\
			\end{array}
		\end{equation}
		and the transformation matrix is \(T=[\delta_4^1\ \delta_4^3\ \delta_4^2\ \delta_4^4]\).
		
		Under the \(z\)-coordinate, the state transition graph is shown by the red parts and black arrows in Fig.~\ref{fig3}~(a), where \(1',2',3',4'\) represent \(\delta_4^1,\delta_4^2,\delta_4^3,\delta_4^4\), respectively. Hence \(z(t+1)=\delta_4[3\ 4\ 2\ 1]z(t)\). From the correspondence in \eqref{eq:xz1}, we obtain \(z^{[2]}(t+1)=\delta_2[1\ 2\ 2\ 1]z(t)\).
	\end{example}
	
	\begin{remark}
		The state transition graph of the \(z^{[1]}\)-subsystem is precisely the quotient graph obtained from BN \eqref{eq:BNAE} by contracting each indistinguishability class into a single vertex.
	\end{remark}
	
	\begin{remark}
		The equivalence classes \(\{1,3\}\) and \(\{2,4\}\) can also be represented by the logical vectors \(\bigl[\begin{smallmatrix}0\\1\end{smallmatrix}\bigr]\) and \(\bigl[\begin{smallmatrix}1\\0\end{smallmatrix}\bigr]\), respectively. Moreover, the values of \(z^{[2]}\) can be labelled by different logical vectors. Different choices yield different transformation matrices \(T\). This is analogous to linear systems, where fixing a basis for the unobservable subspace still allows different orderings of the basis vectors, each yielding a different transformation matrix \(T\).
	\end{remark}
	
	In the linear system \eqref{eq:linear_decomp}, the observable and unobservable parts can always be separated, whereas a BN cannot always admit such a decomposition. This is due to the product structure \(z=z^{[1]}\ltimes z^{[2]}\). For BN \eqref{eq:BNAE}, if \(n_2=0\), then all equivalence classes induced by indistinguishability are singletons and the BN is observable; if \(n_2=n\), then all states belong to a single equivalence class. By Theorem~\ref{thm1}, observability decomposition is realizable in both cases. In the sequel, we focus on the case \(0<n_2<n\), and call such a BN \emph{partially observable}.
	
	From Theorem~\ref{thm1}, a BN admits an observability decomposition exactly when all equivalence classes induced by indistinguishability have the same cardinality. For a partially observable BN, this condition may fail. To obtain a decomposition in such a case, we propose the following extended observability decomposition method for BNs.
	
	\vspace{1ex}
	\noindent\textbf{Extended Observability Decomposition Method for BNs.}
	Let the classes of BN \eqref{eq:BNAE} be \(C_1,\ldots,C_r\), with \(|C_i|=m_i\). Suppose the \(m_i\) are not all equal, and assume without loss of generality that
	\[
	\max_{i\in[1,r]}|C_i|=|C_1|=m_1.
	\]
	For \(i=1,\ldots,r\), set \(k_i=m_1-m_i\). If \(k_i>0\), add \(k_i\) virtual states \(x_{i,1},\ldots,x_{i,k_i}\) to \(C_i\) and let
	\[
	\tilde C_i=C_i\cup\{x_{i,1},\ldots,x_{i,k_i}\};
	\]
	if \(k_i=0\), set \(\tilde C_i=C_i\). For each virtual state \(x_{i,j}\), choose a successor \(\tilde{x}_{i,j}\in\mathcal{N}(C_i)\) and set
	\[
	\mathcal{N}(x_{i,j})=\{\tilde{x}_{i,j}\}, \quad Hx_{i,j}=H(C_i),
	\]
	where \(H(C_i)\) denotes the common output of the states in \(C_i\). By Proposition~\ref{prop1}, \(\mathcal{N}(C_i)\subseteq C_k\) for some \(k\), so \(\tilde{x}_{i,j}\in C_k\). Since all states in \(C_k\) are indistinguishable, the output trajectory from \(\tilde{x}_{i,j}\) coincides with that from the successor of any state in \(C_i\). Together with \(Hx_{i,j}=H(C_i)\), this implies that \(x_{i,j}\) is indistinguishable from every state in \(C_i\). Different classes remain distinguishable: each virtual state is indistinguishable from the states in its own class, while the original classes are pairwise distinguishable, and the transitions among the original states are unchanged. Hence \(\tilde C_1,\ldots,\tilde C_r\) are exactly the indistinguishability classes of the extended logical system, and all of them have cardinality \(m_1\).
	
	The extended system may have a state set whose cardinality is not a power of two. This does not affect the decomposition, since Theorem~\ref{thm1} and its proof only use the logical matrix representation and permutation matrices, which are defined for any finite state set. Since all \(\tilde C_i\) have equal cardinality, Theorem~\ref{thm1} implies that the extended logical system admits an observability decomposition. The state transition graph of the original BN is a subgraph of that of the extended logical system, which implies that the original BN can be embedded into the extended logical system. We call this \emph{the extended decomposition}.
	
	\begin{remark}
		The virtual states have only outgoing edges and no incoming edges, and no initial state of the original BN is a virtual state. Therefore, adding the virtual states does not affect the dynamical evolution of the original BN.
	\end{remark}
	
	\begin{theorem}
		In the sense of extended decomposition, every finite BN admits an observability decomposition.
		\label{thm2}
	\end{theorem}
	
	\begin{example}
		Consider the BN in Example~\ref{eg1}, whose equivalence classes have cardinalities \(1,2,1\). By Theorem~\ref{thm1}, it does not admit an observability decomposition. We now apply the extended decomposition method.
		
		Keep \(\{x_2,x_4\}\) unchanged, add a virtual state \(x_5\) to \(\{x_1\}\), and add a virtual state \(x_6\) to \(\{x_3\}\). We require \(x_5\) to be indistinguishable from \(x_1\), and \(x_6\) to be indistinguishable from \(x_3\). It suffices to set
		\begin{align*}
			\mathcal{N}(x_5) &= \mathcal{N}(x_1)=\{x_3\}, \quad Hx_5 =Hx_1=y^1,\\
			\mathcal{N}(x_6) &= \mathcal{N}(x_3)=\{x_4\}, \quad Hx_6 =Hx_3=y^1.
		\end{align*}
		The extended system is shown in Fig.~\ref{fig4}~(a), where the gray parts are the newly added virtual states.
		
		\begin{figure}[!t]
			\centering
			\begin{minipage}[c]{0.48\textwidth}
				\centering
				\includegraphics[width=5.77cm,height=3.6cm]{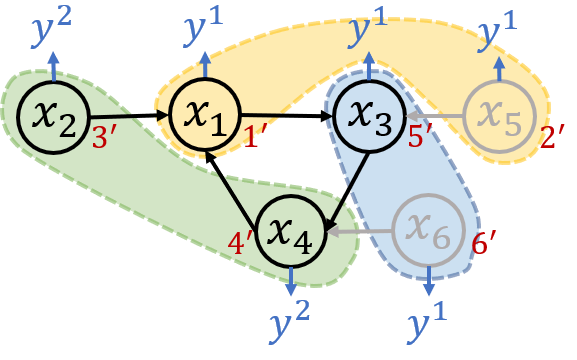}
			\end{minipage}
			\hspace{0.02\textwidth}
			\begin{minipage}[c]{0.48\textwidth}
				\centering
				\includegraphics[width=2.48cm,height=3.5cm]{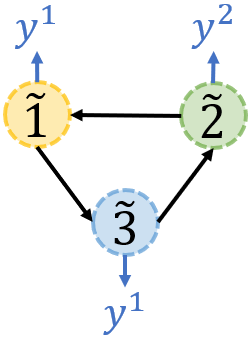}
			\end{minipage}
			\caption{(a) Extended state transition graph of the BN in Example~\ref{eg1}; (b) contracted state transition graph of the extended system \eqref{eq:bne}.}
			\label{fig4}
		\end{figure}
		
		In Fig.~\ref{fig1}, let \(x_i\) represent \(\delta_4^i\), \(i=1,2,3,4\), and \(y^j\) represent \(\delta_2^j\), \(j=1,2\). The original BN is
		\begin{equation*}
			\begin{aligned}
				x(t+1) &= \delta_4[3\ 1\ 4\ 1]x(t),\\
				y(t) &= \delta_2[1\ 2\ 1\ 2]x(t).
			\end{aligned}
		\end{equation*}
		
		In Fig.~\ref{fig4}~(a), let \(x_i\) represent \(\delta_6^i\), \(i=1,\ldots,6\), and \(y^j\) represent \(\delta_2^j\), \(j=1,2\). The extended logical system is
		\begin{equation}
			\begin{aligned}
				\tilde{x}(t+1) &= \delta_6[3\ 1\ 4\ 1\ 3\ 4]\tilde{x}(t), \\
				y(t) &= \delta_2[1\ 2\ 1\ 2\ 1\ 1]\tilde{x}(t).
				\label{eq:bne}
			\end{aligned}
		\end{equation}
		This is a six-state mixed-valued logical network, while Fig.~\ref{fig1} is a part of Fig.~\ref{fig4}~(a). From Fig.~\ref{fig4}~(a), the equivalence classes of \eqref{eq:bne} are \(\{x_1,x_5\}\), \(\{x_2,x_4\}\), and \(\{x_3,x_6\}\), as shown by the yellow, green and blue shaded regions in Fig.~\ref{fig4}~(a). Let \(\delta_3^1,\delta_3^2,\delta_3^3\) label these three classes, respectively. The contracted graph is shown in Fig.~\ref{fig4}~(b), where \(\tilde{i}\) represents \(\delta_3^i\), \(i = 1, 2, 3\).
		
		Let \(\{x_1,x_2,x_3\}\) and \(\{x_5,x_4,x_6\}\) be represented by \(\delta_2^1\) and \(\delta_2^2\), respectively. The coordinate correspondence is
		\begin{equation*}
			\renewcommand{\arraystretch}{1.5}
			\setlength{\arraycolsep}{4pt}
			\begin{array}{ccccccc}
				\tilde{x} & 1 & 2 & 3 & 4 & 5 & 6\\
				z^{[1]} & \delta_3^1 & \delta_3^2 & \delta_3^3 & \delta_3^2 & \delta_3^1 & \delta_3^3\\
				z^{[2]} & \delta_2^1 & \delta_2^1 & \delta_2^1 & \delta_2^2 & \delta_2^2 & \delta_2^2\\
				z & 1' & 3' & 5' & 4' & 2' & 6'
			\end{array}
		\end{equation*}
		and the transformation matrix is
		\[
		T=[\delta_6^1\ \delta_6^3\ \delta_6^5\ \delta_6^4\ \delta_6^2\ \delta_6^6].
		\]
		
		The state transitions under the new coordinate \(z\) are shown in red in Fig.~\ref{fig4}~(a), where \(i'\) represents \(\delta_4^i\), \(i \in [1, 4]\). The state transition equation under the \(z\)-coordinate is
		\[
		z(t+1)=\delta_6[5\ 5\ 1\ 1\ 4\ 4]z(t).
		\]
		Hence the observability decomposition of \eqref{eq:bne} is
		\begin{equation*}
			\begin{aligned}
				z^{[1]}(t+1) &= \delta_3[3\ 1\ 2]z^{[1]}(t),\\
				z^{[2]}(t+1) &= \delta_2[1\ 1\ 1\ 1\ 2\ 2]z(t),\\
				y(t) &= \delta_2[1\ 2\ 1]z^{[1]}(t).
			\end{aligned}
		\end{equation*}
	\end{example}
	
	\begin{remark}
		The observability matrix of BN \eqref{eq:BNAE} is \(\bigl[H^T\ (HL)^T\ \cdots\ (HL^{N-1})^T\bigr]^T\), see \cite{cheng2010analysis}. Since a virtual state inherits the output trajectory of its equivalence class, its column in the observability matrix of the extended system is identical to that of any state in the class. Hence the column comparisons can reuse the information from the original BN, and the additional cost is due only to the added states, not to a complete recomputation.
	\end{remark}
	
	\begin{remark}
		The output decomposition method in \cite{zou2017graph} does not guarantee that the resulting subsystems are observable. For example, suppose the equivalence classes of BN \eqref{eq:BNAE} are \(\{x_1,x_2,x_3,x_4\}\), \(\{x_5,x_6\}\), and \(\{x_7,x_8\}\). The vertex partition produced by the method may be \(\{x_1,x_2\}\), \(\{x_3,x_4\}\), \(\{x_5,x_6\}\), and \(\{x_7,x_8\}\). In this case, the subsystems associated with \(\{x_1,x_2\}\) and \(\{x_3,x_4\}\) are not observable, because the states within each of these two blocks are indistinguishable.
	\end{remark}
	
	\section{Observability Decomposition of BCNs}
	
	A Boolean control network is described by
	\begin{equation}
		\label{eq:bcnle}
		\begin{aligned}
			X(t+1) &= f(X(t),U(t)), \\
			Y(t) &= h(X(t)),
		\end{aligned}
	\end{equation}
	where \(X\in\mathcal{B}^n\), \(U\in\mathcal{B}^m\), and \(Y\in\mathcal{B}^p\) are the state, input, and output variables, and \(f:\mathcal{B}^n\times\mathcal{B}^m\to\mathcal{B}^n\) and \(h:\mathcal{B}^n\to\mathcal{B}^p\) are logical functions.
	
	Let \(x\in\Delta_N\), \(u\in\Delta_M\), and \(y\in\Delta_P\) represent \(X\), \(U\), and \(Y\), respectively, where \(N=2^n\), \(M=2^m\), and \(P=2^p\). Then BCN \eqref{eq:bcnle} can be expressed in algebraic form as \cite{cheng2010analysis}
	\begin{equation}
		\label{eq:bcnae}
		\begin{aligned}
			x(t+1) &= L u(t)x(t), \\
			y(t) &= H x(t),
		\end{aligned}
	\end{equation}
	where \(L\in\mathcal{L}_{N\times NM}\) and \(H\in\mathcal{L}_{P\times N}\). Here and below, the semi-tensor product is omitted.
	
	Write \(L=[L_1\ L_2\ \cdots\ L_M]\) with \(L_s\in\mathcal{L}_{N\times N}\), \(s\in[1,M]\). At time \(t\), taking \(u(t)=\delta_M^s\) gives \(x(t+1)=L_s x(t)\). Thus, for a fixed input sequence, BCN \eqref{eq:bcnae} can be regarded as a switched system
	\[
	x(t+1)=L_{\sigma(t)}x(t),
	\]
	where \(\sigma(t)=s\) when \(u(t)=\delta_M^s\). The state transition graph of \(x(t+1)=L_s x(t)\) is denoted by \(\mathcal{G}_s\) and is called the \(s\)-th state transition graph of BCN \eqref{eq:bcnae}, \(s=1,\ldots,M\).
	
	\begin{definition}[{\cite{fornasini2012observability}}]
		\label{def:2}
		BCN \eqref{eq:bcnae} is said to be observable if, for any two distinct initial states \(x_0=\delta_N^i\) and \(\bar{x}_0=\delta_N^j\) with \(i\neq j\), there exists an input sequence \(\{u(t)\}_{t\ge 0}\) such that the output sequences \(\{y(t)\}_{t\ge 0}\) and \(\{\bar{y}(t)\}_{t\ge 0}\) are distinct. Two states \(x_0\) and \(\bar{x}_0\) are said to be distinguishable if such an input sequence exists; otherwise, they are said to be indistinguishable.
	\end{definition}
	
	The indistinguishability relation introduced by Definition~\ref{def:2} is again an equivalence relation, which partitions the state set of BCN \eqref{eq:bcnae} into equivalence classes; let \(\mathcal{C}_\sim\) denote the set of these classes. Similar to Proposition~\ref{prop1}, we have the following property.
	
	\begin{proposition}
		\label{prop2}
		Given BCN \eqref{eq:bcnae}, let \(C_1,C_2,\ldots,C_r\) be the equivalence classes induced by the indistinguishability relation. Then, for any \(C_i\), there exists an index \(k\) such that
		\[
		\mathcal{N}_s(C_i)\subseteq C_k, \forall s = 1, \ldots, M.
		\]
	\end{proposition}
	
	\begin{proof}
		Suppose, for the sake of contradiction, that there exist \(C_i\), \(i \in [1, r]\), and \(s \in [1, M]\), such that the conclusion does not hold. The following proof is by contradiction, similar to that of Proposition~\ref{prop1}.
	\end{proof}
	
	Following the linear case, we define observability decomposition for BCNs as follows.
	
	\begin{definition}[{\cite{li2022observability}}]
		\label{def:bcn_decomp}
		BCN \eqref{eq:bcnae} admits an observability decomposition if there exists a coordinate transformation \(z=Tx\) with a permutation matrix \(T\) such that, in the \(z\)-coordinate frame, the BCN becomes
		\begin{equation}
			\begin{aligned}
				z^{[1]}(t+1) &= G_1 u(t) z^{[1]}(t), \\
				z^{[2]}(t+1) &= G_2 u(t) z(t), \\
				y(t) &= Q z^{[1]}(t),
			\end{aligned}
			\label{eq:bcndecomp}
		\end{equation}
		where the \(z^{[1]}\)-subsystem is observable, \(z^{[1]}\in\Delta_{N_1}\), \(z^{[2]}\in\Delta_{N_2}\), \(u\in\Delta_M\), \(y\in\Delta_P\), \(z=z^{[1]}\ltimes z^{[2]}\), and \(N_1=2^{n_1}\), \(N_2=2^{n_2}\), \(M=2^m\), \(P=2^p\), \(N_1N_2=N\). The matrices satisfy \(G_1\in\mathcal L_{N_1\times N_1M}\), \(G_2\in\mathcal L_{N_2\times NM}\), and \(Q\in\mathcal L_{P\times N_1}\).
	\end{definition}
	
	Write \(G_1=[G_{11}\ G_{12}\ \cdots\ G_{1M}]\) and \(G_2=[G_{21}\ G_{22}\ \cdots\ G_{2M}]\), where \(G_{1s}\in\mathcal L_{N_1\times N_1}\) and \(G_{2s}\in\mathcal L_{N_2\times N}\) for \(s=1,\ldots,M\). Then \(u(t)\) acts as a switching signal in the first two equations of \eqref{eq:bcndecomp}.
	
	The following lemma relates the matrices of BCN \eqref{eq:bcnae} and those of its decomposition \eqref{eq:bcndecomp}.
	
	\begin{lemma}[{\cite{zou2017graph}}]
		\label{lem:2}
		If BCN \eqref{eq:bcnae} admits an observability decomposition, then \((L,H)\) and \((G_1,G_2,Q)\) satisfy
		\begin{align*}
			&(\bm I_{N_1}\otimes\bm 1_{N_2}^T)TL_sT^T = G_{1s}(\bm I_{N_1}\otimes\bm 1_{N_2}^T), \\
			&(\bm 1_{N_1}^T\otimes\bm I_{N_2})TL_sT^T = G_{2s}, \\
			&(\bm 1_{N_1}^T\otimes\bm I_{N_2})TL(\bm I_M\otimes T^T)=G_2, \\
			&HT^T=Q(\bm I_{N_1}\otimes\bm 1_{N_2}^T).
		\end{align*}
	\end{lemma}
	
	\begin{theorem}
		\label{thm3}
		BCN \eqref{eq:bcnae} admits an observability decomposition if and only if all equivalence classes induced by the indistinguishability relation have the same cardinality.
	\end{theorem}
	
	\begin{proof}
		\emph{[Necessity]}
		The argument is the same as in Theorem~\ref{thm1}, with the input fixed at each step. If BCN \eqref{eq:bcnae} admits a decomposition \eqref{eq:bcndecomp}, then states with the same \(z^{[1]}\)-coordinate are indistinguishable, while states with different \(z^{[1]}\)-coordinates are distinguishable because the \(z^{[1]}\)-subsystem is observable. Hence the indistinguishability classes all have the same cardinality.
		
		\emph{[Sufficiency]}
		Suppose the indistinguishability classes are \(C_1,\ldots,C_{N_1}\), each of cardinality \(N_2\), where \(N_1N_2=N\). We construct \(T,G_1,Q,G_2\) satisfying \eqref{eq:bcndecomp}.
		
		\emph{Construction of \(T\).}
		The matrix \(T\) is constructed as in the proof of Theorem~\ref{thm1}; it is a permutation matrix satisfying \eqref{t} and let \(z = T x\).
		
		\emph{Construction of \(G_1\).}
		For each \(l=1,\ldots,N_1\), by Proposition~\ref{prop2} there exists an index \(k\) such that \(\mathcal{N}_s(C_l)\subseteq C_k\), \(s \in [1, M]\). Define
		\begin{equation}
			Col_l(G_{1s})=\delta_{N_1}^{k},
			\qquad l=1,\ldots,N_1.
			\label{eq:g1s}
		\end{equation}
		The choice of k is independent of the representative chosen in \(C_l\) and subgraph \(\mathcal{G}_s\), since the image of the whole class \(C_l\) under any \(L_s\) lies in \(C_{k}\). If \(x(t)=\delta_N^{a_i^l}\) and \(u(t)=\delta_M^s\), then \(z^{[1]}(t)=\delta_{N_1}^l\) and \(x(t+1)=L_s\delta_N^{a_i^l}\in C_{k}\), so
		\[
		z^{[1]}(t+1)=\delta_{N_1}^k=G_{1s}\delta_{N_1}^l=G_{1s}z^{[1]}(t).
		\]
		Hence \(z^{[1]}(t+1)=G_1u(t)z^{[1]}(t)\).
		
		\emph{Construction of \(Q\).}
		The matrix \(Q\) is constructed as in the proof of Theorem~\ref{thm1}. Since the output equation does not depend on the input, the same argument gives
		\[
		y(t)=Qz^{[1]}(t).
		\]
		
		\emph{Construction of \(G_2\).}
		Let
		\begin{equation}
			G_{2s}=(\bm 1_{N_1}^T\otimes\bm I_{N_2})TL_sT^T,
			\qquad s=1,\ldots,M,
			\label{eq:g2s}
		\end{equation}
		and set \(G_2=[G_{21}\ G_{22}\ \cdots\ G_{2M}]\).
		At time \(t\), let the state of the \(x\)-system be \(\delta_N^a\), the state of the \(z\)-system be \(T\delta_N^a\), and the input be \(u(t)=\delta_M^s\). Then at time \(t+1\), the \(x\)-state is \(L_s\delta_N^a\) and the \(z\)-state is \(TL_s\delta_N^a\). By \eqref{t}, there exist \(l\in[1,N_1]\) and \(i\in[1,N_2]\) such that
		\[
		T(L_s\delta_N^a)=\delta_{N_1}^l\ltimes\delta_{N_2}^i.
		\]
		Therefore,
		\[
		\delta_{N_2}^i
		=(\bm 1_{N_1}^T\otimes\bm I_{N_2})T(L_s\delta_N^a)
		=(\bm 1_{N_1}^T\otimes\bm I_{N_2})TL_sT^T(T\delta_N^a)
		=G_{2s}(T\delta_N^a),
		\]
		which implies \(z^{[2]}(t+1)=G_{2s}z(t)=G_2u(t)z(t)\).
		
		It remains to show that the \(z^{[1]}\)-subsystem is observable. Suppose, for contradiction, that two distinct states \(\delta_{N_1}^l\) and \(\delta_{N_1}^k\) of the \(z^{[1]}\)-subsystem are indistinguishable. Then for every input sequence,
		\[
		QG_{1s_{t-1}}\cdots G_{1s_0}\delta_{N_1}^l
		=QG_{1s_{t-1}}\cdots G_{1s_0}\delta_{N_1}^k,\qquad t\ge 1.
		\]
		Take any \(x_l\in C_l\) and \(x_k\in C_k\). Their \(z^{[1]}\)-coordinates are \(\delta_{N_1}^l\) and \(\delta_{N_1}^k\), and their output trajectories in the original BCN coincide for every input sequence, contradicting the fact that \(C_l\) and \(C_k\) are different indistinguishability classes. Hence the \(z^{[1]}\)-subsystem is observable.
		
		Therefore, the constructed \(T,G_1,Q,G_2\) satisfy \eqref{eq:bcndecomp} and the \(z^{[1]}\)-subsystem is observable, so BCN \eqref{eq:bcnae} admits an observability decomposition.
		\hfill\(\square\)
	\end{proof}
	
	\vspace{1ex}
	\noindent\textbf{Extended Observability Decomposition Method for BCNs.}
	Let the equivalence classes of BCN \eqref{eq:bcnae} be \(C_1,\ldots,C_r\), with \(|C_i|=m_i\). Suppose the \(m_i\) are not all equal, and assume without loss of generality that
	\[
	\max_{i\in[1,r]}|C_i|=|C_1|=m_1.
	\]
	For \(i=1,\ldots,r\), set \(k_i=m_1-m_i\). If \(k_i>0\), add \(k_i\) virtual states \(x_{i,1},\ldots,x_{i,k_i}\) to \(C_i\) and let
	\[
	\tilde C_i=C_i\cup\{x_{i,1},\ldots,x_{i,k_i}\};
	\]
	if \(k_i=0\), set \(\tilde C_i=C_i\). For each \(s=1,\ldots,M\) and each virtual state \(x_{i,j}\), choose a successor \(\tilde{x}_{i,j,s}\in\mathcal{N}_s(C_i)\) and set
	\[
	\mathcal{N}_s(x_{i,j})=\{\tilde{x}_{i,j,s}\},\quad H(x_{i,j})=H(C_i),
	\]
	where \(H(C_i)\) denotes the common output of the states in \(C_i\). By Proposition~\ref{prop2}, for all \(s\) there exists an index \(k\) such that \(\mathcal{N}_s(C_i)\subseteq C_{k}\), so \(\tilde{x}_{i,j,s}\in C_{k}\). Since all states in \(C_{k}\) are indistinguishable, the output trajectory from \(\tilde{x}_{i,j,s}\) under any input coincides with that from the successor of any state in \(C_i\). Together with \(Hx_{i,j}=H(C_i)\), this implies that \(x_{i,j}\) is indistinguishable from every state in \(C_i\). Different classes remain distinguishable: each virtual state is indistinguishable from the states in its own class, while the original classes are pairwise distinguishable, and the transitions among the original states are unchanged. Hence \(\tilde C_1,\ldots,\tilde C_r\) are exactly the indistinguishability classes of the extended logical system, and all of them have cardinality \(m_1\).
	
	The extended system may have a state set whose cardinality is not a power of two, so it is a mixed-valued logical system. This does not affect the decomposition, since Theorem~\ref{thm3} and its proof only use the logical matrix representation and permutation matrices, which are defined for any finite state set. Since all \(\tilde C_i\) have equal cardinality, Theorem~\ref{thm3} implies that the extended logical system admits an observability decomposition. Moreover, every state transition graph of the original BCN is a part of the corresponding state transition subgraph of the extended logical system.
	
	\begin{theorem}
		In the sense of extended decomposition, every finite BCN admits an observability decomposition.
		\label{thm4}
	\end{theorem}
	
	The equivalence classes in the preceding examples can be read directly from the graphs. For a general BCN, they have to be computed from the observability matrices \cite{cheng2009controllability,fornasini2012observability}. Define the first, second, third, \(\ldots\) layer observation matrices by
	\begin{equation*}
		\mathcal T_1=H,\qquad
		\mathcal T_2=
		\begin{bmatrix}
			HL_1\\
			\vdots\\
			HL_M
		\end{bmatrix},\qquad
		\mathcal T_3=
		\begin{bmatrix}
			\mathcal T_2L_1\\
			\vdots\\
			\mathcal T_2L_M
		\end{bmatrix},
		\quad\ldots
	\end{equation*}
	In general,
	\[
	\mathcal T_{k+1}=
	\begin{bmatrix}
		\mathcal T_kL_1\\
		\vdots\\
		\mathcal T_kL_M
	\end{bmatrix},
	\qquad k\ge 1.
	\]
	Equivalently, in the STP form, \(\mathcal T_{k+1}=\mathcal T_k\ltimes[L_1^T\ \cdots\ L_M^T]^T\), where the product is understood blockwise. The matrix \(\mathcal T_k\) collects the output blocks \(HL_{s_{k-1}}\cdots L_{s_1}\) over all input sequences \(\delta_M^{s_1},\ldots,\delta_M^{s_{k-1}}\) of length \(k-1\), and its columns correspond to the \(N\) initial states.
	
	By comparing the columns of \(\mathcal T_1\), the one-step equivalence classes \(\mathcal C_{\sim,1}\) are obtained. For the columns of \(\mathcal T_2\) corresponding to the states in each class of \(\mathcal C_{\sim,1}\), pairwise comparison gives the two-step equivalence classes \(\mathcal C_{\sim,2}\). This is because \(\mathcal C_{\sim,2}\) refines \(\mathcal C_{\sim,1}\). The process continues: suppose \(\mathcal C_{\sim,k}=\{C_{k,1},\ldots,C_{k,r}\}\). Compute \(\mathcal T_{k+1}\), and compare, within each \(C_{k,i}\), the columns of \(\mathcal T_{k+1}\) corresponding to the states in \(C_{k,i}\). This yields \(\mathcal C_{\sim,k+1}\). When \(\mathcal C_{\sim,k}=\mathcal C_{\sim,k+1}\), the partition has stabilized and equals the overall equivalence classes \(\mathcal C_\sim\).
	
	We present an algorithm below to realize the observability decomposition of BCN \eqref{eq:bcnae}.
	
	\begin{algorithm}
		\footnotesize
		\caption{Implementation of (extended) observability decomposition for BCNs.}
		\label{alg}
		\begin{algorithmic}[1]
			\REQUIRE \(L=[L_1\ L_2\ \cdots\ L_M]\in\mathcal L_{N\times NM}\), \(L_s\in\mathcal L_{N\times N}\), \(H\in\mathcal L_{P\times N}\);
			\ENSURE A coordinate transformation \(T\) and matrices \(G_1,G_2,Q\) of the (extended) decomposition;
			
			\STATE \(\mathcal T_1\Leftarrow H\);
			\STATE Compute the one-step equivalence classes \(\mathcal C_{\sim,1}\) by pairwise comparison of the columns of \(\mathcal T_1\);
			
			\FOR{\(k=2\) to \(N\)}
			\STATE \(\mathcal T_k\Leftarrow
			\mathcal T_{k-1}
			\begin{bmatrix}
				L_1\\
				\vdots\\
				L_M
			\end{bmatrix}\);
			\STATE Compute the \(k\)-step equivalence classes \(\mathcal C_{\sim,k}\) by pairwise comparison of the columns of \(\mathcal T_k\) within each class of \(\mathcal C_{\sim,k-1}\);
			\IF{\(\mathcal C_{\sim,k}=\mathcal C_{\sim,k-1}\)}
			\STATE \(\mathcal C_\sim\Leftarrow\mathcal C_{\sim,k}\);
			\STATE \textbf{break}
			\ENDIF
			\ENDFOR
			
			\IF{all classes in \(\mathcal C_\sim\) have cardinality \(1\)}
			\STATE BCN \eqref{eq:bcnae} is observable;
			\STATE \textbf{return} ``observable''
			\ELSE
			\IF{all classes in \(\mathcal C_\sim\) have the same cardinality \(N_2>1\)}
			\STATE Construct \(T,G_1,G_2,Q\) as in the proof of Theorem~\ref{thm3};
			\STATE \textbf{return} \(T,G_1,G_2,Q\)
			\ELSE
			\STATE Add virtual states to each class so that all classes have equal cardinality, following the extended decomposition method;
			\STATE Construct \(T,G_1,G_2,Q\) for the extended logical system as in the equal-cardinality case;
			\STATE \textbf{return} the extended logical system and \(T,G_1,G_2,Q\)
			\ENDIF
			\ENDIF
		\end{algorithmic}
	\end{algorithm}
	
	\begin{remark}
		Considering the dimensions of \(L\) and \(H\), Algorithm~\ref{alg} runs in exponential time. However, the partition \(\mathcal C_{\sim,k-1}\) is reused, so pairwise comparisons are performed only within each class of \(\mathcal C_{\sim,k-1}\), which reduces the number of column comparisons.
	\end{remark}
	
	The following example illustrates the observability decomposition and the extended decomposition of BCNs.
	
	\begin{example}
		We consider a simplified model of the \textit{Escherichia coli} lac operon \cite{veliz2011boolean}.
		\begin{equation}
			\begin{aligned}
				X_1(t+1) &= X_2(t)\vee X_3(t), \\
				X_2(t+1) &= U_1(t)\wedge X_1(t), \label{lacl}\\
				X_3(t+1) &= U_1(t)\vee (U_2(t)\wedge X_1(t)), \\
				Y(t) &= X_1(t),
			\end{aligned}
		\end{equation}
		where \(X_1,X_2,X_3,U_1,U_2\) represent the mRNA, the LacZ polypeptide, internal lactose, external lactose, and external lactose metabolite, respectively.
		
		The algebraic form of \eqref{lacl} is
		\begin{equation}
			\begin{aligned}
				x(t+1) &= \delta_8[1\ 1\ 1\ 5\ 3\ 3\ 3\ 7 \mid 1\ 1\ 1\ 5\ 3\ 3\ 3\ 7 \mid 3\ 3\ 3\ 7\ 4\ 4\ 4\ 8\mid 4\ 4\ 4\ 8\ 4\ 4\ 4\ 8]u(t)x(t), \\
				y(t) &= \delta_2[1\ 1\ 1\ 1\ 2\ 2\ 2\ 2]x(t).
				\label{laca}
			\end{aligned}
		\end{equation}
		Its state transition subgraphs \(\mathcal G_1,\mathcal G_2,\mathcal G_3,\mathcal G_4\) are shown in Fig.~\ref{fig51}, where \(i\) represents \(\delta_8^i\), \(i=1,\ldots,8\), and \(y^j\) represents \(\delta_2^j\), \(j=1,2\).
		
		\begin{figure}[!t]
			\centering
			\begin{minipage}[c]{0.24\textwidth}
				\centering
				\includegraphics[width=3.5cm,height=3.2cm]{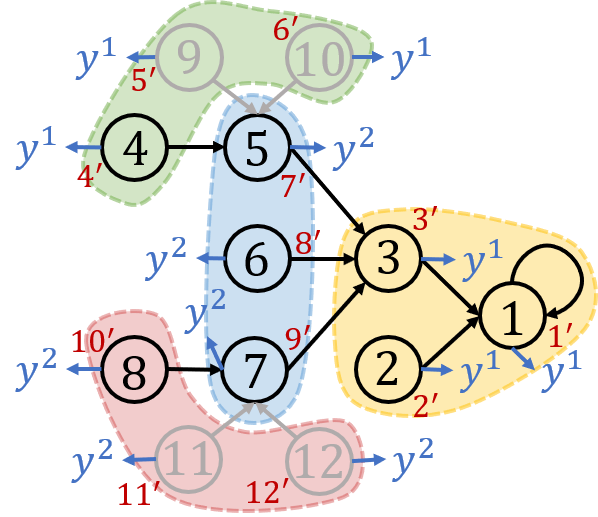}
			\end{minipage}
			\hspace{0.00\textwidth}
			\begin{minipage}[c]{0.24\textwidth}
				\centering
				\includegraphics[width=3.5cm,height=3.2cm]{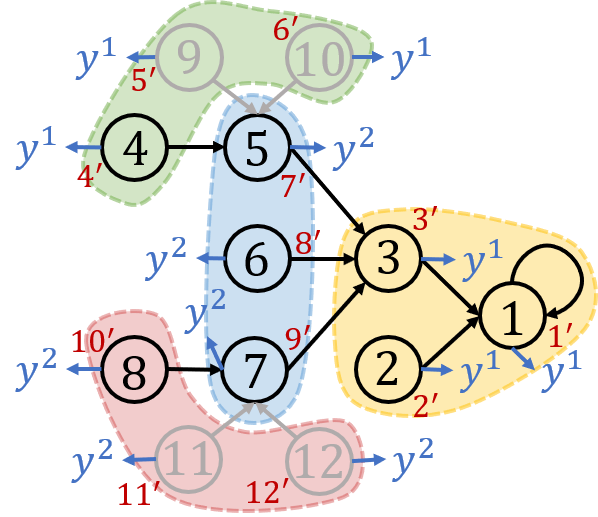}
			\end{minipage}
			\hspace{0.00\textwidth}
			\begin{minipage}[c]{0.24\textwidth}
				\centering
				\includegraphics[width=3.8cm,height=3.1cm]{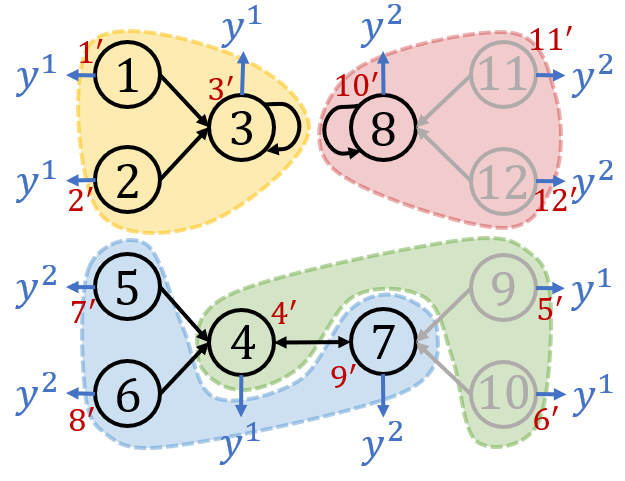}
			\end{minipage}
			\hspace{0.00\textwidth}
			\begin{minipage}[c]{0.24\textwidth}
				\centering
				\includegraphics[width=3.5cm,height=3.4cm]{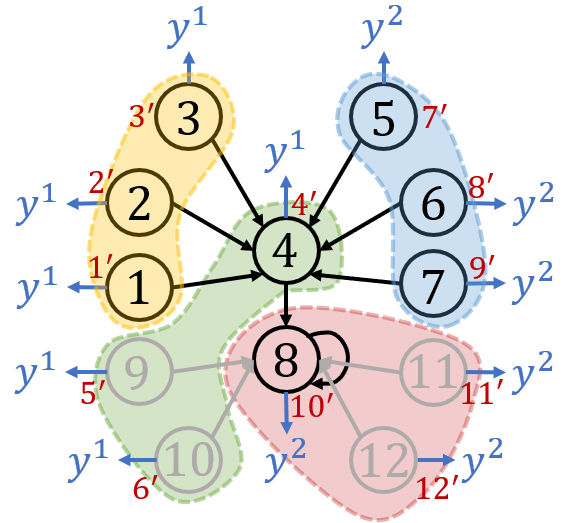}
			\end{minipage}
			\caption{State transition subgraphs of \eqref{laca} and \eqref{es}; gray parts are virtual states. (a) \(\mathcal{G}_1\); (b) \(\mathcal{G}_2\); (c) \(\mathcal{G}_3\); (d) \(\mathcal{G}_4\).}
			\label{fig51}
		\end{figure}
		
		The first observation matrix is \(\mathcal T_1=\delta_2[1\ 1\ 1\ 1\ 2\ 2\ 2\ 2]\), so the one-step equivalence classes are
		\[
		\mathcal C_{\sim,1}=\bigl\{\{1,2,3,4\},\{5,6,7,8\}\bigr\}.
		\]
		The second observation matrix is
		\[
		\mathcal T_2=
		\begin{bmatrix}
			\delta_2[1\ 1\ 1\ 2\ 1\ 1\ 1\ 2]\\
			\delta_2[1\ 1\ 1\ 2\ 1\ 1\ 1\ 2]\\
			\delta_2[1\ 1\ 1\ 2\ 1\ 1\ 1\ 2]\\
			\delta_2[1\ 1\ 1\ 2\ 1\ 1\ 1\ 2]
		\end{bmatrix},
		\]
		so the two-step equivalence classes are
		\[
		\mathcal C_{\sim,2}=\bigl\{\{1,2,3\},\{4\},\{5,6,7\},\{8\}\bigr\}.
		\]
		Computing \(\mathcal T_2L_1,\mathcal T_2L_2,\mathcal T_2L_3,\mathcal T_2L_4\) yields the four blocks shown below
		\begin{gather*}
			\resizebox{\linewidth}{!}{$
				\begin{bmatrix} \delta_2[1\ 1\ 1\ 1\ 1\ 1\ 1\ 1] \\ \delta_2[1\ 1\ 1\ 1\ 1\ 1\ 1\ 1] \\ \delta_2[1\ 1\ 1\ 1\ 1\ 1\ 1\ 1] \\ \delta_2[1\ 1\ 1\ 1\ 1\ 1\ 1\ 1] \end{bmatrix}, \quad
				\begin{bmatrix} \delta_2[1\ 1\ 1\ 1\ 1\ 1\ 1\ 1] \\ \delta_2[1\ 1\ 1\ 1\ 1\ 1\ 1\ 1] \\ \delta_2[1\ 1\ 1\ 1\ 1\ 1\ 1\ 1] \\ \delta_2[1\ 1\ 1\ 1\ 1\ 1\ 1\ 1] \end{bmatrix}, \quad
				\begin{bmatrix} \delta_2[1\ 1\ 1\ 1\ 2\ 2\ 2\ 2] \\ \delta_2[1\ 1\ 1\ 1\ 2\ 2\ 2\ 2] \\ \delta_2[1\ 1\ 1\ 1\ 2\ 2\ 2\ 2] \\ \delta_2[1\ 1\ 1\ 1\ 2\ 2\ 2\ 2] \end{bmatrix}, \quad
				\begin{bmatrix} \delta_2[2\ 2\ 2\ 2\ 2\ 2\ 2\ 2] \\ \delta_2[2\ 2\ 2\ 2\ 2\ 2\ 2\ 2] \\ \delta_2[2\ 2\ 2\ 2\ 2\ 2\ 2\ 2] \\ \delta_2[2\ 2\ 2\ 2\ 2\ 2\ 2\ 2] \end{bmatrix},
				$}
		\end{gather*}
		and hence
		\[
		\mathcal T_3=[(\mathcal T_2L_1)^T\ (\mathcal T_2L_2)^T\ (\mathcal T_2L_3)^T\ (\mathcal T_2L_4)^T]^T.
		\]
		The three-step equivalence classes are
		\[
		\mathcal C_{\sim,3}=\bigl\{\{1,2,3\},\{4\},\{5,6,7\},\{8\}\bigr\}.
		\]
		Since \(\mathcal C_{\sim,3}=\mathcal C_{\sim,2}\), the equivalence classes of \eqref{laca} are
		\[
		\mathcal C_\sim=\bigl\{\{x_1,x_2,x_3\},\{x_4\},\{x_5,x_6,x_7\},\{x_8\}\bigr\},
		\]
		where \(x_i\) represents \(\delta_8^i\), \(i\in[1,8]\).
		
		Since the cardinalities \(3,1,3,1\) are not all equal, Theorem~\ref{thm3} does not apply; we now apply the extended decomposition method. We add virtual states \(x_9,x_{10}\) to \(\{x_4\}\), and virtual states \(x_{11},x_{12}\) to \(\{x_8\}\). To make \(\{x_4,x_9,x_{10}\}\) and \(\{x_8,x_{11},x_{12}\}\) equivalence classes, it suffices to choose, for \(i=1,\ldots,4\), successors \(\tilde{x}_{4,i}\in\mathcal N_i(x_4)\) and \(\tilde{x}_{8,i}\in\mathcal N_i(x_8)\), and set
		\begin{align*}
			\mathcal N_i(x_9) &= \mathcal N_i(x_{10})=\{\tilde{x}_{4,i}\}, \quad Hx_9=Hx_{10}=Hx_4=y^1,\\
			\mathcal N_i(x_{11}) &= \mathcal N_i(x_{12})=\{\tilde{x}_{8,i}\}, \quad Hx_{11}=Hx_{12}=Hx_8=y^2.
		\end{align*}
		
		\begin{figure}[!t]
			\centering
			\begin{minipage}[c]{0.24\textwidth}
				\centering
				\includegraphics[width=3.5cm,height=3.9cm]{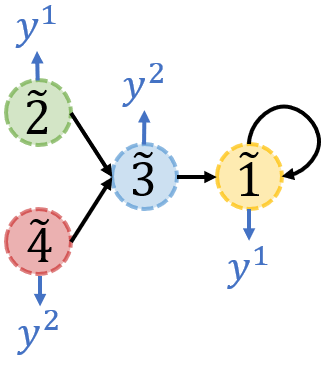}
			\end{minipage}
			\hspace{0.00\textwidth}
			\begin{minipage}[c]{0.24\textwidth}
				\centering
				\includegraphics[width=3.5cm,height=3.9cm]{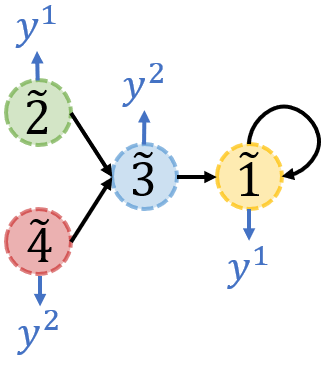}
			\end{minipage}
			\hspace{0.00\textwidth}
			\begin{minipage}[c]{0.24\textwidth}
				\centering
				\includegraphics[width=3.5cm,height=3.95cm]{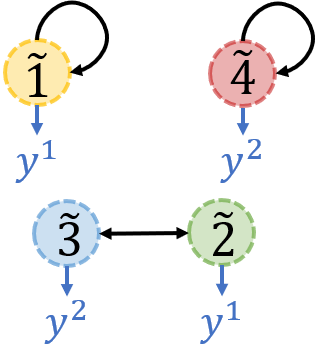}
			\end{minipage}
			\hspace{0.00\textwidth}
			\begin{minipage}[c]{0.24\textwidth}
				\centering
				\includegraphics[width=3.5cm,height=4.1cm]{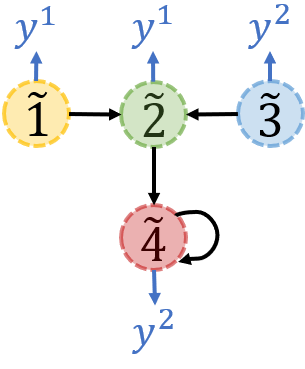}
			\end{minipage}
			\caption{Contracted state transition subgraphs of the extended system \eqref{es}: (a) \(\mathcal{G}_1\); (b) \(\mathcal{G}_2\); (c) \(\mathcal{G}_3\); (d) \(\mathcal{G}_4\).}
			\label{fig52}
		\end{figure}
		
		The gray parts in Fig.~\ref{fig51} are the virtual states. The extended system is
		\begin{equation}
			\begin{aligned}
				\tilde{x}(t+1) &= \delta_{12}[1\ 1\ 1\ 5\ 3\ 3\ 3\ 7\ 5\ 5\ 7\ 7\mid 1\ 1\ 1\ 5\ 3\ 3\ 3\ 7\ 5\ 5\ 7\ 7\mid \\
				&\qquad\quad 3\ 3\ 3\ 7\ 4\ 4\ 4\ 8\ 7\ 7\ 8\ 8\mid 4\ 4\ 4\ 8\ 4\ 4\ 4\ 8\ 8\ 8\ 8\ 8]u(t)\tilde{x}(t), \\
				y(t) &= \delta_2[1\ 1\ 1\ 1\ 2\ 2\ 2\ 2\ 1\ 1\ 2\ 2]\tilde{x}(t).
				\label{es}
			\end{aligned}
		\end{equation}
		The equivalence classes of \eqref{es} are \(\{x_1,x_2,x_3\}\), \(\{x_4,x_9,x_{10}\}\), \(\{x_5,x_6,x_7\}\), \(\{x_8,x_{11},x_{12}\}\), as shown by the yellow, green, blue and red shaded regions in Fig.~\ref{fig51}. Let \(\delta_4^1,\delta_4^2,\delta_4^3,\delta_4^4\) in \(z^{[1]}\) label these four classes, respectively. The contracted graph of \eqref{es} is shown in Fig.~\ref{fig52}, where \(\tilde i\) represents \(\delta_4^i\), \(i=1,2,3,4\). Let \(\{x_1,x_4,x_5,x_8\}\), \(\{x_2,x_9,x_6,x_{11}\}\), \(\{x_3,x_{10},x_7,x_{12}\}\) represent \(\delta_3^1,\delta_3^2,\delta_3^3\) in \(z^{[2]}\), respectively. Then the coordinate correspondence is
		\begin{equation*}
			\renewcommand{\arraystretch}{1.5}
			\setlength{\arraycolsep}{4pt}
			\begin{array}{ccccccccccccc}
				\tilde{x} & 1 & 2 & 3 & 4 & 5 & 6 & 7 & 8 & 9 & 10 & 11 & 12\\
				z^{[1]} & \delta_4^1 & \delta_4^1 & \delta_4^1 & \delta_4^2 & \delta_4^3 & \delta_4^3 & \delta_4^3 & \delta_4^4 & \delta_4^2 & \delta_4^2 & \delta_4^4 & \delta_4^4\\
				z^{[2]} & \delta_3^1 & \delta_3^2 & \delta_3^3 & \delta_3^1 & \delta_3^1 & \delta_3^2 & \delta_3^3 & \delta_3^1 & \delta_3^2 & \delta_3^3 & \delta_3^2 & \delta_3^3\\
				z & 1' & 2' & 3' & 4' & 7' & 8' & 9' & 10' & 5' & 6' & 11' & 12'
			\end{array}
		\end{equation*}
		and the transformation matrix is
		\[
		T=[\delta_{12}^1\ \delta_{12}^2\ \delta_{12}^3\ \delta_{12}^4\ \delta_{12}^7\ \delta_{12}^8\ \delta_{12}^9\ \delta_{12}^{10}\ \delta_{12}^5\ \delta_{12}^6\ \delta_{12}^{11}\ \delta_{12}^{12}].
		\]
		
		The state transition graphs under the \(z\)-coordinate are shown in red in Fig.~\ref{fig51}, where \(i'\) represents \(\delta_{12}^i\) for \(i\in[1,12]\). The \(z\)-state evolution is
		\begin{equation*}
			\begin{aligned}
				z(t+1) &= \delta_{12}[1\ 1\ 1\ 7\ 7\ 7\ 3\ 3\ 3\ 9\ 9\ 9\mid 1\ 1\ 1\ 7\ 7\ 7\ 3\ 3\ 3\ 9\ 9\ 9\mid\\
				&\qquad\quad 3\ 3\ 3\ 9\ 9\ 9\ 4\ 4\ 4\ 10\ 10\ 10\mid 4\ 4\ 4\ 10\ 10\ 10\ 4\ 4\ 4\ 10\ 10\ 10]u(t)z(t).
			\end{aligned}
		\end{equation*}
		Hence the extended observability decomposition of \eqref{laca} is
		\begin{equation*}
			\begin{aligned}
				z^{[1]}(t+1) &= \delta_4[1\ 3\ 1\ 3\mid 1\ 3\ 1\ 3\mid 1\ 3\ 2\ 4\mid 2\ 4\ 2\ 4]u(t)z^{[1]}(t),\\
				z^{[2]}(t+1) &= \delta_3[1\ 1\ 1\ 1\ 1\ 1\ 3\ 3\ 3\ 3\ 3\ 3\mid
				1\ 1\ 1\ 1\ 1\ 1\ 3\ 3\ 3\ 3\ 3\ 3 \mid\\
				&\quad\quad\;\; 3\ 3\ 3\ 3\ 3\ 3\ 1\ 1\ 1\ 1\ 1\ 1\mid
				1\ 1\ 1\ 1\ 1\ 1\ 1\ 1\ 1\ 1\ 1\ 1]u(t)z(t),\\
				y(t) &= \delta_2[1\ 1\ 2\ 2]z^{[1]}(t).
			\end{aligned}
		\end{equation*}
		The contracted graph in Fig.~\ref{fig52} also shows directly that the \(z^{[1]}\)-subsystem is observable.
	\end{example}
	
	\section{Conclusion}
	
	This paper studied observability decomposition of BNs and BCNs via indistinguishability classes. For a BN, we proved that an observability decomposition exists if and only if all such equivalence classes have the same cardinality, and we gave an explicit construction of the transformation matrix and the subsystem matrices. We also showed that the state transition graph of the observable subsystem is the quotient graph obtained by contracting each equivalence class into a single vertex. For BNs whose equivalence classes have unequal cardinalities, we proposed an extended decomposition by adding virtual states, and showed that every finite BN can be embedded into an extended logical system that admits an observability decomposition. These results were extended to BCNs, and an algorithm together with an illustrative example was provided.

	\Acknowledgements{This work was supported by the National Natural Science Foundation of China (Grant No. 62373193).}


	\DeclarationofAI{During the preparation of this manuscript, the authors used DeepSeek-V4.1-Flash-1017 for text translation. The authors have manually verified the AI-generated content and take full responsibility for the accuracy, completeness, and originality of the entire manuscript.}


	\bibliographystyle{plain}
	\bibliography{myrefers}
	
	
	
\end{document}